\documentstyle[12pt]{article}
\title{On the Distribution of Zeros of a Ruelle Zeta-Function\thanks{
AMS Subject Classification 58F20, 30D05, 30D20. Keywords:
Iteration, Zeta function, Ruelle operator, Entire function,
Subharmonic function}}
\author{A. Eremenko\thanks{Supported by NSF grant DMS-9101798}, G. Levin
and M. Sodin}
\date{December 27, 1992}

\newtheorem{theorem}{Theorem}
\newtheorem{lemma}{Lemma}

\newtheorem{proposition}{Proposition}
\begin{document}
\maketitle
\begin{abstract}
We study the limit distribution of zeros of a Ruelle $\zeta$-function
for the dynamical system
$z\mapsto z^2+c$ when $c$ is real and $c\rightarrow -2-0$
and apply the results to the correlation functions of this dynamical
system.
\end{abstract}
	
Consider the dynamical system defined by the complex polynomial map
$f_c:z\mapsto z^2+c$, where $c<-2$. We use the notions and results of 
the iteration theory of rational functions (see for example \cite{EL}). 
Denote by $f_c^{*n}$ the $n-$th
iterate of the function $f_c$. The Julia set $J(f_c)$ is a
Cantor set on the real line. So in particular all finite
periodic points are real. 
This system is expanding (hyperbolic) on its Julia set.
When $c=-2$ the Julia set is the segment $[-2,2]$ and
the map $P=f_{-2}$ is not expanding anymore. We have the 
conjugation
\begin{equation}
\label{conjugation}
P\circ\phi=\phi\circ Q,
\end{equation}
where $
\phi:[0,1]\rightarrow [-2,2],\;t\mapsto 2\cos\pi t$ and
$$
Q=\left\{\begin{array}{ll}
t\mapsto 2t, & 0\leq t\leq 1/2,\\
t\mapsto 2-2t, & 1/2\leq t\leq 1.
\end{array}
\right.
$$
 Remark that the chaotic dynamic of $P$ on $[-2,2]$
was investigated by J. von Neuman and S. Ulam on one of the first computers.

We are going to study the dynamics of $f_c,\;c<-2$ when $c\rightarrow -2$
and then compare it with the behavior of the limit system $P$.
The chaotic dynamics of $f_c$ has to be described in probabilistic
terms. This can be done by introducing an appropriate
invariant probability measure $\sigma_c$ on the Julia set.
We will show that the rate of asymptotic decrease of correlation functions
of the system $(f_c,\nu_c)$ changes dramatically when we pass to the limit
system as $c\rightarrow -2$.

Our tool is the Thermodynamic Formalism \cite{Ru0,Ru,Ru1,Ru2}.
Let us introduce the main objects of this theory in our particular case.
%VSTAVKA
Consider the Fr\'{e}chet space $C^{\infty}(U)$ of infinitely differentiable
functions defined in some real neighborhood $U$ of the Julia set, such that
$f_c^{-1}(U)\subset U$ and $U$ does not contain the critical point of $f_c$.
%First
We define the Ruelle operator $L_c$ acting on $C^{\infty}(U)$
%on the space of holomorphic functions in a neighborhood of
%the Julia set:
by the formula
$$L_cg(x)=
\sum_{\{y:f_c(y)=x\}}\frac{g(y)}{[f'_c(y)]^2}.$$
The weight $(f'_c)^{-2}$ is strictly
positive on $U$.
%KONEC VSTAVKI
According to Ruelle's extension of the Perron-Frobenius theorem 
$L_c$ has a simple maximal positive eigenvalue $\lambda^{-1}_0(c)$ such that
the moduli of all other eigenvalues are strictly less
then $|\lambda^{-1}_0(c)|$.
Let $h_c$ and $\nu_c$ denote the eigenvectors of $L_c$ and the conjugate
operator $L^*_c$ respectively, corresponding to the eigenvalue
$\lambda^{-1}_0(c)$ ($h_c$ is a positive continuous function and $\nu_c$ is a
Borel measure).
Then $\sigma_c=h_c\nu_c$ is an $f_c$-invariant ergodic probability measure
on the Julia set, called
%VSTAVKA
``The Gibbs state, corresponding to
the weight $(f'_c)^{-2}$''. 
%The following particular form of Ruelle's zeta-function is connected
%to the operator $L_c$:
%$$\zeta_c(\lambda)=\exp{ \left( \sum_{m=1}^{\infty}\frac{\lambda^m}{m}
%\sum_{x\in {\rm Fix}(f^m_c)}\frac{1}{(f_c^{*m})'(x)} \right)},$$ 
%where ${\rm Fix}(f^m_c)$ is the set of fixed points of $f^{*m}_c$.
%VSTAVKA
%(We chose the weight $\phi=({f'_c})^{-1}$ in the definition of Ruelle
% $\zeta$-function. See section 8 of
%\cite{Ru1} and formula (3.3) with $\sigma=\infty$ in \cite{LSY}.)
%VSTAVKA
The operator $L_c$ can be also considered on the space $A$
of functions analytic in a complex neighborhood of the Julia set.
%ISMENENIE.CRAVNI C TEM,CHTO ECT U TEBIA
Namely, for every complex neighborhood $W$ of the Julia set such that
$U\subset W$, $f_c^{-1}(W)\subset W$ and $W$ does not contain the critical point of
$f_c$, consider the Banach space $A(W)$ of functions analytic in $W$ with the supremum norm. Then $A$ is the union of all such $A(W)$.
%KONEZ ISMENENII
As the weight $(f'_c)^{-2}$ is analytic, the spectrum and
eigenfunctions of $L_c$ in $A$ are the same as in $C^{\infty}(U)$
(see \cite{Ru1}, Corollary 3.3(i)). This fact allows us to use
the explicit expressions for eigenfunctions found in \cite{LSY}
with the help of complex analysis.
%Pereneceno c predidushei crtanizi.GENA
The following particular form of Ruelle's zeta-function is connected
to the operator $L_c$:
$$\zeta_c(\lambda)=\exp{ \left( \sum_{m=1}^{\infty}\frac{\lambda^m}{m}
\sum_{x\in {\rm Fix}(f^m_c)}\frac{1}{(f_c^{*m})'(x)} \right)},$$ 
where ${\rm Fix}(f^{*m}_c)$ is the set of fixed points of $f^{*m}_c$.
%VSTAVKA
(We chose the weight $\phi=({f'_c})^{-1}$ in the definition of Ruelle
 $\zeta$-function. See section 8 of
\cite{Ru1} and formula (3.3) with $\sigma=\infty$ in \cite{LSY}.)
%VSTAVKA
The function $\zeta_c$ can be expressed in terms of generalized
Fredholm determinants (\cite[Corollary 8.1]{Ru1}).
In our particular case it coincides with the Fredholm deteminant
$D_c$ of $L_c$ \cite{LSY}; this is an entire function of order zero and
its zeros are reciprocal to the eigenvalues of $L_c$. There is an
explicit formula found in \cite{LSY} (see also \cite{LSY2}):
$$\zeta_c(\lambda)=D_c(\lambda)=
1+\sum_{n=1}^{\infty}\frac{\lambda^n}{2^nf_c(0)\ldots 
f^{*n}_c(0)}.$$
In Appendix we will give a short direct proof of the fact that the
eigenvalues of $L_c$ are reciprocal to the zeros of $D_c$.

%GENA
{\sl Remark.} Let us consider another extension of the operator $L_c:C^\infty (U)\to
C^\infty (U)$ to the Fr\'{e}chet space $C^\infty (W)$ of the $C^\infty$-functions of two real variables $u$ and $v$, $u+iv\in W$, given by the formula
$$L_c^{\bf R^2} g(x)=
\sum_{\{y:f_c(y)=x\}}\frac{g(y)}{|f'_c(y)|^2}.$$
(Note that $|f'_c(x)|^2$ is the Jacobian of the map $f'_c:\bf R^2\to \bf R^2$
at the point $x$.) Then the eigenvalues and eigenfunctions of $L_c^{\bf R^2}$ coincide with those for $L_c$. Really, every eigenfunction of $L_c^{\bf R^2}$ restricted to $U=W\bigcap \bf R$ is an eigenfunction of $L_c$. Conversely, the eigenfunctions of $L_c$ are analytic and, hence, belong to $C^\infty (W)$. In particular, $(\lambda_0(c))^{-1}$ is the leading eigenvalue of the operator $L_c^{\bf R^2}$ and the value $\log \lambda_0(c)$ is the so-called ``escape rate'' \cite{KT}. 
%END

One of the reasons why the study of eigenvalues of $L_c$ is important is
their connection to correlation functions.
For any two continuous $A$ and $B$ on the Julia set
define the correlation function $\rho_{c,A,B}$ by
$$\rho_{c,A,B}(m)=\sigma_c(A(f^{*m}_c).B)-\sigma_c(A).\sigma_c(B),$$
where $\sigma(A)=\int A\,d\sigma$. Let
$$S_{c,A,B}(z)=\sum_{m=0}^{\infty}\rho_{c,A,B}(m)z^m$$
be the corresponding generating function. {\em If $A$ and $B$ are
infinitely differentiable 
on the Julia set then $S_{c,A,B}$ is meromorphic
in ${\bf C}$ and its poles can be located only at the points $\lambda
\lambda^{-1}_0$, where $\lambda^{-1}$ runs over the 
eigenvalues of $L_c$ other than}
$\lambda_0$ \cite[Proposition 5.3]{Ru1}. 

{\bf 1}. First we investigate the limit distribution of eigenvalues of $L_c$ or,
which is equivalent, zeros of $D_c$. 
The following facts about distribution of zeros of $D_c$ were
established in \cite{L}. 
For all $c<-2$ the zeros with moduli greater than $1000$ are negative,
and simple. There exists a constant $c_0=-2.85\ldots$ such that
for $c\leq c_0$ all zeros of $D_c$ are real. If $c<-2$
is close to $-2$ then there are 
non-real zeros and their number tends to infinity as $c$ tends to $-2$.

To study the asymptotic distribution of complex zeros we introduce
the probability measures $\mu_c$ which charge equally every zero
whose modulus is less than $1000$.
\begin{theorem}.
The measures $\mu_c$ tend weakly to the uniform distribution
on the circle $\{\lambda:|\lambda|=4\}$.
\end{theorem}
{\sl Remarks.} Notice that $4$ is the radius of convergence of the
series $D_{-2}=(4-2\lambda)/(4-\lambda)$.
Our proof is also applicable to the family of 
entire functions
$$H_a(z)=\sum_{n=0}^{\infty}\frac{z^n}{a^{2^n-1}},\quad a>1,$$
whose distribution of zeros was studied by G. H. Hardy \cite{H}.
He proved that for fixed $a$ all zeros with moduli greater than
$r_0(a)$ are negative. (In fact $r_0(a)$ can
be replaced by an absolute constant
\cite{L}).
Our argument shows that
the limit distribution of zeros of $H_a$ when $a\rightarrow 1$
is the uniform distribution on the circle $\{z:|z|=1\}$.
Theorem 1 should be compared with the following theorem of
Jentzsch and Szeg\H{o}:
{\em the limit distribution of zeros of partial sums
of a power series $\sum a_kz^k$ is the uniform distribution
on $\{z:|z|=1\}$, provided that $|a_k|^{1/k}\rightarrow 1$.
}  Our proof is based on the same idea as Beurling's 
proof of the Jentzsch-Szeg\H{o} theorem \cite{B}.

{\em Proof.} We assume that $-3<c<-2$. 
It is convenient to introduce the variable $z=\lambda/2$ and
set $F_c(z)=D_c(2z)$ and $r_n(c)=f_c^{*n}(0)$. Thus
$$F_c(z)=1+\sum_{n=1}^{\infty}\frac{z^n}{r_1(c)\ldots r_n(c)}$$
and 
\begin{equation}
\label{rec}
r_{n+1}(c)=r_n^2(c)+c,\qquad r_1(c)=c,\quad c<-2.
\end{equation}
It is easy to see that 
all $r_n$, except $r_1$, are positive, the sequence $(r_n)$ 
is increasing and
$r_{n+1}(c)/r_n(c)\rightarrow \infty,\;n\rightarrow  \infty,
\; c<-2.$ Denote by $k=k(c)$ the smallest natural $k$ such that
\begin{equation}
\label{>36}
\frac{r_{k+1}(c)}{r_k(c)}\geq 36.
\end{equation}
It was proved in \cite{L} that the number of zeros of $F_c$ in 
any fixed disk $\{z:|z|<R\},\;R>1000$ is asymptotically equivalent to
$k(c)$ when $c\rightarrow -2.$ This fact also follows from the 
estimates below (formula (7) plus Rouch\'{e} theorem).
\begin{lemma}.
If $k=k(c)$ is as defined above, then

{\bf (i)}$\;\;\; 36\leq |r_k(c)|\leq 1521=39^2$.

{\bf (ii)}$\;\; k(c)\sim (\log |c+2|^{-1})/{\log 4},\quad c\rightarrow -2.$

{\bf (iii)} $(1/k(c))\log |r_1(c)\ldots r_{k(c)}(c)|\rightarrow \log2,
\quad c\rightarrow -2.$
\end{lemma}
{\em Proof.} {\bf (i)} From (\ref{>36}) we conclude that $k=k(c)>1$. 
If $r_k(c)<36$ then by (\ref{rec}) 
$r_{k+1}(c)/r_k(c)=r_k(c)+c/r_k(c)<r_k(c)<36$,
which contradicts the definition of $k$. 
This proves the left inequality in (i).
Now assume that $|r_k|>39^2$. Then in view of (\ref{rec})
we have $|r_{k-1}|>39$ and we obtain $|r_k|=|r_{k-1}|^2+c>|r_{k-1}|^2-3$ and
$|r_k|/|r_{k-1}|>36$, which contradicts the definition of $k$. This proves
the right inequality in (i).

{\bf (ii)} Set $c=-2-t,\;t>0$. An easy induction gives
\begin{equation}
\label{lower}
|r_n(c)|\geq 2+(4^{n-1}-1)t,\quad n=1, 2, \ldots .
\end{equation}
To prove an inequality in the opposite direction we remark that
$r_{n+1}(c)=[r_n(c)]^2-2-t\leq [r_n(c)]^2-2=P(r_n(c))$, so 
$$r_n(c)\leq P^{*(n-1)}
(r_1(c))\leq P^{*(n-1)}(2+t).$$
Using the semiconjugacy
$$2\cosh 2z=[2\cosh z]^2-2=P(2\cosh z),$$
(it is more convenient to use $\cosh$ rather then $\cos$ here)
we obtain $r_n(c)\leq 2\cosh (2^{n-1}y),$ where $y$ is the smallest 
positive solution of the equation
such that $2\cosh y=2+t$. There exists an 
absolute constant $C_0=30$ such that $2\cosh x \leq C_0 x^2+2$ whenever
$2\cosh x\leq 1521,\;x\in {\bf R}$. Thus we obtain
\begin{equation}
\label{upper}
r_n(c)\leq 2+4^{n-1}C_0 t,\quad n=1, 2, \ldots, k(c).
\end{equation}
The statement (ii) follows from (\ref{lower}) and (\ref{upper}).

{\bf (iii)} 
 From (ii) follows
\begin{equation}
\label{t}
t\leq C_1 4^{-k}.
\end{equation}
In view of (\ref{lower}), (\ref{upper}) (\ref{t}) we have
$$\left| \left(\frac{1}{k}\sum_{n=1}^k\log |r_n(c)|\right)
-\log 2\right| \leq
\frac{1}{k}\sum_{n=1}^k\log (1+4^{n-1}C_0t)\leq$$
$$\leq\frac{1}{k}\sum_{n=1}^kC_0 C_1 4^{n-k}\leq
\frac{1}{k}\sum_{n=0}^{\infty}C_0 C_1 
4^{-n}\rightarrow 0,\quad k\rightarrow\infty.$$
This finishes the proof of the Lemma 1.

Denote by $A(t_1,t_2)$ the annulus $\{z:t_1<|z|<t_2\}$ and set
$A(c)=A(4r_k(c),9r_k(c)),$ where $k=k(c)$. 
Put $M_c(z)=z^k/(r_1(c)\ldots r_k(c)).$
If $z\in A(c)$ we have
\begin{equation}
\label{est}
\left| 1-\frac{F_c(z)}{M_c(z)}\right| \leq 
\sum_{j=1}^k \frac{r_k\ldots r_{k-j+1}}{
|z|^j}+\sum_{j=1}^{\infty}
\frac{|z|^j}{r_{k+1}\ldots r_{k+j}}\leq$$
$$\leq\sum_{j=1}^{\infty} 4^{-j} + \sum_{j=1}^{\infty} 4^{-j}=\frac{2}{3}.
\end{equation}
Thus if we denote $u_c(z)=(k(c))^{-1}\log |F_c(z)|$ 
then by (iii) of the Lemma 1
\begin{equation}
\label{approx}
u_c(z)=(k(c)^{-1})\log |M_c(z)|+o(1)=\log |z/2|+o(1),\quad
c\rightarrow -2,
\end{equation}
uniformly when $z\in A(c)$.
We are going to prove that 
\begin{equation}
\label{convergence}
u_c(z)\rightarrow \log^+|z/2|,\quad |z|\leq 324,
\end{equation}
where the convergence holds in $L^1$ with respect to the 
Lebesgue measure (area)
in $\{z:|z|\leq 324\}$.

 From the definition of $A(c)$ and Lemma 1, (i) follows that
$A_c\subset
A(144,13689)$. So from any sequence $c_m\rightarrow -2,\;c_m<-2$
we can chose a subsequence (which we again denote 
by $c_m$ such that the annuli 
$A(c_m)$ contain
a fixed annulus $A(q_1,q_2),\; q_1<q_2,\; q_2>324.$
Then in view of (\ref{approx}) we have 
\begin{equation}
\label{fixapprox}
u_{c_m}(z)\rightarrow \log |z/2|\quad {\rm uniformly\; in}\; 
{\bar A(q_1,q_2)}.
\end{equation}
Furthermore we have
\begin{equation}
\label{approxunit}
u_{c_m}(z)\rightarrow 0,\quad |z|<2 
\end{equation}
(convergence in $L^1$ on compacts in $\{z:|z|<2\}$),
because $F_c(z)\rightarrow F_{-2}(z)=1-z/(2-z),\;c\rightarrow -2$
uniformly on compacts in $\{z:|z|<2\}$.
Now we use the following fact (see for example 
\cite{Hor}, Theorem 4.1.9): {\em if a sequence of subharmonic
functions $u_m$ is bounded from above on $\{z:|z|=R\}$ and their values
at the point $0$ are bounded from below then there is a subsequence
which converges in $L^1$ on every compact in $\{z:|z|<R\}$ to a
subharmonic function $u$.} Applying this statement to our functions
$u_{c_m}$ and $R=q_2$, we obtain a subsequence (which we again denote by
$u_{c_m}$) which converges to a subharmonic function $u$. This function
$u$ has the properties:
\begin{equation}
\label{zero}
u(z)=0,\quad |z|<2
\end{equation}
and
\begin{equation}
\label{log}
u(z)=\log |z/2|,\quad q_1<|z|<q_2,
\end{equation}
which follow from (\ref{approxunit}) and (\ref{fixapprox}) respectively.
Remark that $u(z)\leq 0,\;|z|=2$. This follows from (\ref{zero})
and the following theorem of M. Brelot \cite{Brelot}: {\em
if $u$ is a subharmonic function and $u(z_0)=a$ then for every
$\epsilon>0$ there exists
a sequence of circles centered at $z_0$ and radii tending to zero
such that $u(z)\geq a-\epsilon$ on these circles.} 
(It follows from the upper semi-continuity of $u$ that 
$u(z)\geq 0,\;|z|=2$, but we do not need this.)
Now $\log |z/2|$ is a harmonic majorant of $u$ in the annulus $A(2,q_2)$,
but $u(z)=\log |z/2|$ at some points in this annulus, for example for
$|z|=q_1$. It follows from the Maximum Principle that 
$u(z)=\log^+ |z/2|,\;|z|<q_2$.

Thus we have proved that from {\em every} sequence $u_{c_m}$ we can select
a subsequence tending to $\log^+ |z/2|$. This means that (\ref{convergence})
is true. In fact our proof shows that $u_c$ converge
to $\log^+|z/2|$ in $L^1$ on every compact in the plane.
Now we conclude from the general results on convergence of
subharmonic functions \cite{A,AB,Hor} that the Riesz measures $\mu_c$ of $u_c$ 
converge weakly to the Riesz measure of $u$, which is the uniform
measure on the circle $|\lambda|=2|z|=4$. This proves the theorem.
\vspace{.2in}

{\bf 2}. Now we consider the application of Theorem 1 to the
dynamical system $(f_c,\sigma_c)$ where $\sigma_c$
is the Gibbs state defined in the introduction. We have
$$\zeta_c(\lambda)\rightarrow 1-\frac{\lambda}{4-\lambda},\quad
c\rightarrow -2,$$
uniformly on compacts in $\{\lambda:|\lambda|<4\}$. So $\lambda_0(c)
\rightarrow 2$ and 
$$\inf\{\lambda:\zeta_c(\lambda)=0,\;
\lambda\neq \lambda_0\}\rightarrow 4,\quad c\rightarrow -2.$$
Thus by Theorem 1 and by Ruelle's theorem mentioned in introduction we have the following
asympotic behavior of correlation functions:
$$\limsup_{m\rightarrow\infty}|\rho_{c,A,B}(m)|^{1/m}=r(c),$$
$$\mbox{where $r(c)\rightarrow 1/2$ as $c\rightarrow -2.$}$$
We want to compare this result with the behavior of the limiting dynamical system
when $c\rightarrow -2$. First we have to understand what the limit
invariant measure is. Recall the conjugation 
(\ref{conjugation}).
The Lebesgue measure $l_1$ on $[0,1]$ 
is invariant with respect to $Q$ thus its image $\sigma_{-2}=\phi_*l_1$ is
invariant with respect to $P=f_{-2}$. The measure $\sigma_{-2}$ is absolutely
continuous with the density
$$\frac{1}{\pi \sqrt{4-x^2}}$$
on the interval $[-2,2]$.
\begin{proposition}.
$\sigma_c\rightarrow\sigma_{-2}$ weakly as $c\rightarrow -2$.
\end{proposition}
{\em Proof}. We will use the explicit expressions for the eigenfunction
$h_c$ of $L_c$ and for the Cauchy transform 
$$H_c(z)=\int\frac{d\nu_c(x)}{x-z}$$
of the eigenmeasure $\nu_c$ of $L^*_c$, corresponding to the greatest
eigenvalue $\lambda^{-1}_0$ (see \cite{SY,LSY}). Using the notation
$r_n(c)=f_c^{*n}(0)$ we have
$$h_c(x)=\sum_{n=0}^{\infty}\frac{\lambda_0^n(c)}{2^nr_1(c)\ldots r_n(c)
[r_{n+1}(c)-x]}$$
and
$$H_c(z)=
\sum_{n=0}^{\infty}
\frac{\lambda_0^n(c)}{2^nzf_c(z)\ldots f_c^{*n}(z)}.$$
The function $z\mapsto H_c(z)$ is holomorphic in
the complement
of the Julia set $J(f_c)$.
We have
$$h_c(x)\rightarrow -\left(\frac{1}{2+x}+\frac{1}{2-x}\right),\quad c\rightarrow -2$$
in ${\bf \bar C}\backslash ((-\infty,-2]\cup [2,\infty))$
and
$$H_c(z)\rightarrow H_{-2}(z)=\sum_{n=0}^{\infty}
\frac{1}{zP(z)\ldots P^{*n}(z)},\quad c\rightarrow -2$$
in ${\bf \bar C}\backslash [-2,2]$.   

Consider the measure $\nu_{-2}$ on $[-2,2]$ with the density $\sqrt{4-x^2}$.
We claim that $H_{-2}(z)$ is proportional
to the Caushy transform of $\nu_{-2}$.
This follows from the fact that they both
satisfy the same functional equation
$$H(z)-\frac{H(P(z))}{z}=\frac{{\rm const}}{z},\quad z\in
{\bf \bar C}\backslash [-2,2].$$
Now Proposition 1 follows from the identity
$$\left(\frac{1}{2+x}+\frac{1}{2-x}\right)\sqrt{4-x^2}=
\frac{4}{\sqrt{4-x^2}}.$$

 So the dynamical system $(P,\sigma_{-2})$ is the limit of
$(f_c,\sigma_c)$ when $c\rightarrow -2$. We will show that the
asymptotic behavior
of correlations changes drastically when we pass to the limit as 
$c\rightarrow-2.$ 
\begin{proposition}.
Let $A$ and $B$ be holomorphic functions on $[-2,2]$.
Then there exists a constant $a=a(A,B)>1$ such that
$$\rho_{-2,A,B}(m)\sim a^{-2^m},\quad m\rightarrow\infty.$$
\end{proposition}
{\em Proof}. In view of Cauchy formula is enough to prove the
proposition for the set of functions
$$A_z(x)=\frac{1}{z-x},\quad x\in [-2,2],\quad z\in {\bf \bar C}
\backslash [-2,2].$$
After the pullback to the segment $[0,1]$ via the conjugation
(\ref{conjugation}) we have to consider the correlations
$$\rho_{A,B}(m)=l_1(A(Q^m).B)-l_1(A).l_1(B)$$
with $A$ and $B$ of the form
$$\frac{1}{z-2\cos\pi t}.$$
If we introduce the operator
\begin{equation}
\label{defg}
G:\;g(t)\mapsto \frac{1}{2}\sum_{y:Q(y)=t}g(y)=\frac{1}{2}
\left(g(t/2)+g(1-t/2)\right)
\end{equation}
then
\begin{equation}
\label{G}
\rho_{A,B}(m)=l_1(A.G^m(B))-l_1(A).l_1(B).
\end{equation}
Now we notice that
$$G\left(\frac{1}{z-2\cos\pi t}\right)=\frac{P'(z)}{2(P(z)-2\cos\pi t)},$$
which implies
\begin{equation}
\label{analytic}
G^m\left(\frac{1}{z-2\cos\pi t}\right)=\frac{(P^{*m})'(z)}{2^m(P^{*m}(z)-
2\cos\pi t)}=
S(z)+\frac{\cos\pi t+o(1)}{2^{m-1}(P^{*m}(z))^2},
\end{equation}
where $S$ is a function depending only on $z$.
Combining (\ref{G}) and (\ref{analytic}) we get the statement of
Proposition 2.

{\sl Remark.}
The analyticity assumption in Proposition 2 is crucial.  
Indeed consider the operator $G$ defined in (\ref{defg})
in the space of infinitely differentiable functions on
$[0,1]$. Its eigenvalues are $4^{-m}\;m=0,1,2\ldots,$ and to each
eigenvalue $4^{-m}$ corresponds one (up to
a constant multiple) eigenfunction $p_m$ which is a polynomial
of degree $2m$.
Now if $A$ and $B$ belong to the subspace of $L^2([0,1],l_1)$
generated by $\{p_m:m=0,1,2,\ldots\}$ then we have
$$\rho_{A,B}(m)\sim {\rm const}.4^{-km},\quad m\rightarrow\infty,$$
where ${\rm const\neq 0}$ and $k$ depend on $A$ and $B$.

{\bf Appendix.} Here we indicate a direct proof of the fact that
the eigenvalues of $L_c$ are reciprocal to the zeros of $D_c,\;c<-2$
(see also \cite{LSY2}).
Let us look at the eigenvalues of the adjoint operator $L^*_c$.
The dual space $A^*$ is the space of functions $g$ analytic in the complement
of the Julia set $J(f_c)$ and equal to zero at infinity. To every
such function corresponds a linear functional given by
$$h\mapsto\frac{1}{2\pi i}\int gh,$$
where the integral is taken along some countur surrounding $J(f_c)$.
Now a change of the variable in this integral shows that
$\lambda^{-1}$ is an eigenvalue iff for every function $h$
holomorphic in a neighborhood of $J(f_c)$
$$\int \left(g-\lambda \frac{g\circ f_c}{{f_c}'}\right) h=0.$$
Thus $w=g-\lambda g\circ f_c/{f_c}'$ is holomorphic on $J_c$.
It is also holomorphic in ${\bf \bar C}\backslash (J(f_c)\cup \{0\})$ 
because $f^\prime_c(z)=2z$. We conclude that $w(z)={\rm const}/z$
and after the normalization of $g$ we get the functional equation
$$g(z)=\frac{\lambda}{2z}g(f_c(z))+\frac{1}{z},$$
from which follows that
$$g(z)=\frac{1}{z}+\sum_{n=1}^{\infty}\frac{\lambda^n}{2^nzf_c(z)
\ldots f^{*n}_c(z)}.$$
Now $g$ is holomorphic at $0$ so the residue of the series in
the right side should vanish that is
$$D_c(\lambda)=1+\sum_{n=1}^{\infty}\frac{\lambda^n}{2^n f_c(0)
\ldots f^{*n}(0)}=0.$$

We thank the referee for his valuable comments.

{\em Purdue University, West Lafayette IN 47907}

{\em Institute of Mathematics, Hebrew University, Jerusalem, 91904}

{\em Institute of Low Temperature Physics \& Engineering,
 
Kharkov, 310164, Ukraine}

\end{document}